# A mean field Jacobi process for modeling sustainable tourism

Hidekazu Yoshioka [a, *]

[a] Japan Advanced Institute of Science and Technology, 1-1 Asahidai, Nomi, Ishikawa 923-1292, Japan
[*] Corresponding author: yoshih@jaist.ac.jp, ORCID: 0000-0002-5293-3246

**Abstract**
A mean field Jacobi process governing the dynamics of the travel demand of agents is formulated and its application to sustainable tourism is investigated both mathematically and computationally. The bounded nature of the Jacobi diffusion process enables the categorization of tourism state with sustainable tourism state corresponding to an internal solution and overtourism state to a boundary solution. A stochastic control framework is introduced to design sustainable tourism under uncertainty, incorporating model distortion conditions owing to misspecification. The control problem is reduced to solving the optimality system of a stationary mean field game whose closed-form solution is derived under certain conditions. The optimality system can be computed numerically using the finite difference method under more general conditions. We present demonstrative examples of the mean field Jacobi process for different parameter values, illustrating both sustainable tourism and overtourism cases. Our findings suggest that the sustainable tourism state cannot be realized if the fluctuation or model misspecification is large.

**Statements and Declarations**
**Acknowledgements:** This study was supported by the Japanese Society for the Promotion of Science (22K14441 and 22H02456) and Nippon Life Insurance Foundation (Environmental Research Grant for Young Researchers in 2024, No. 24).
**Funding:** This study was supported by the Japanese Society for the Promotion of Science (Grant numbers 22K14441 and 22H02456) and Nippon Life Insurance Foundation (Environmental Research Grant for Young Researchers in 2024, No. 24).
**Data availability statement:** The data will be made available upon reasonable request to the corresponding author.
**Competing interests:** The authors have no competing interests to declare that are relevant to the content of this article.
**Declaration of generative artificial intelligence (AI) in scientific writing:** The authors did not use generative AI technology to prepare this manuscript.
**Ethics approval and consent to participate:** Not applicable.

## 1. Introduction

### 1.1 Problem background

Sustainability is a key concept in the modern world. It refers to the harmonious coexistence between the environment and humanity so that both can be balanced globally in the long run (Vanelli and Kobiyama, 2021; Liu et al., 2023; Han et al., 2024). In this study, we focus on tourism where people and nature interact. Mass tourism has damaged the atmosphere, oceans, freshwater (Buckley, 2012), and world heritage sites. Visitors exceeding the tourism-carrying capacity of a site destroy its essential qualities (Li et al., 2021; Santos and Brilha, 2023). Tourism that is expected to be well sustained in the future, namely sustainable tourism, has recently been an important topic in the fields of environmental, economic, and industrial research. Lee and Chang (2015) proposed a threshold number of tourists required to protect the water environment and ecosystem of a tourist spot; this threshold can be considered as a capacity limit to determine the sustainability of tourism.

Mathematical modeling can be used to analyze sustainable tourism. Chenavaz et al. (2022) reported that a dynamic pricing policy that accounts for the damage caused by visitors and entrance fees for heritage conservation could provide a stronger foundation for sustainable tourism policies. Brida et al. (2013) formulated an economic model and investigated its control problem using a dynamic programming principle to determine the conditions under which economic growth, tourism development, and environmental sustainability are balanced. By operating a dynamic model, Inchausti-Sintes (2023) found that tourism serves as an essential contributor to the preservation of the environment provided that the environment is integrated into the economic system. Antoci et al. (2022) proposed a predator–prey model to explain overtourism. Metilelu et al. (2022) established a system of nonlinear ordinary differential equations to determine the stability of tourism-polluting coastal sites. Multiagent evolutionary game models have been investigated, focusing on the collective behavior emerging among local stakeholders and tourists (Chica et al., 2022; Chica et al., 2023).

Many mathematical models have been proposed for analyzing sustainable tourism; however, those based on the assumption that economic and social dynamics arise from interactions among agents have not been studied extensively except for few studies (Chica et al., 2022; Chica et al., 2023). The mean field game (MFG) which can deal with complex dynamic decision-making among agents through a statistical mechanical principle, meets this requirement (Lasry and Lions, 2007). The statistical mechanical nature of the MFG enables the reduction of infinite-dimensional agent dynamics into a low-dimensional system inheriting essential agent interactions. In this regard, Bagagiolo et al. (2019) followed by Andria et al. (2023) applied MFGs to the sustainable management of tourist flows on a graph network. As a related study, Centorrino et al. (2021) computationally analyzed a crowd management problem in museums as destinations of tourists. However, other research examples are lacking; hence, there is a need to deepen knowledge regarding mathematical modeling in the context of sustainable tourism. Moreover, to the best of our knowledge, an analytically tractable model, namely a model with a closed-form solution, has not been studied. If this model type exists, the conditions under which tourism sustainability is achieved can be determined. These backgrounds motivated the present study, as explained in the following subsection.

**1.2 Aim and contribution**

This study aims to formulate an analytically tractable MFG model for designing sustainable tourism. The contributions of this research are the following.

Our first contribution is the formulation of the tourism dynamics of agents using stochastic differential equations (SDEs) (Capasso and Bakstein, 2021). We focus on the Jacobi process, also called Wright–Fisher and Pearson diffusion, as an analytically tractable SDE whose statistical characteristics, including stability conditions, moments, and stationary probability density functions (PDFs), are found explicitly (Chapter 6 in Alfonsi, 2015; Ditlevsen et al., 2020). The solution to the Jacobi process is bounded in a closed interval. We exploit this property to define the boundary states as no- and overtourism states, while the internal states are sustainable tourism states. In this way, we can efficiently categorize tourism states by using a single SDE. We investigate the classical Jacobi process and its controlled version using an MFG. Because Jacobi processes arise in many engineering problems, such as neuronal dynamics (D'Onofrio et al., 2019; D'Onofrio et al., 2024), energy transition dynamics (Aïd et al., 2021), financial engineering (Tong and Liu, 2022), wind power forecasting (Caballero et al., 2021), and generative machine learning (Avdeyev et al., 2023), this study's outcomes would contribute not only to modeling the sustainable tourism but also to these diverse research areas.

Our second contribution is the proposal of an MFG that admits a closed-form solution under certain conditions, thus providing an analytical characterization of the success or failure of sustainable tourism. Indeed, closed-form solutions to MFGs play a pivotal role in understanding social and economic dynamics such as the linear-quadratic game (Bardi, 2012), commodity market dynamics (Brown and Ambrose, 2024), electricity market trading (Coskun and Korn, 2024), and graphon games (Foguen-Tchuendom et al., 2024). Solving an MFG based on the SDEs of agent dynamics reduces to the identification of proper solutions to a Hamilton–Jacobi–Bellman (HJB) equation that governs the optimal control and worst-case distortion and a Fokker–Planck (FP) equation that governs the PDF of the controlled system (Lasry and Lions, 2007).

Our MFG is unique because it accounts for model distortion, namely misspecification, based on Knight's formalism. According to this formalism, distortion is introduced through a Radon–Nikodym derivative between a benchmark and distorted model (Hansen and Sargent, 2001). An advantage of this formalism is that it harmonizes with the stochastic control theory in which misspecification can be accounted for by an objective function without assuming a particular form of uncertainty distribution. This control-based formalism has recently been applied to assessing uncertainties in climate change scenarios (Xepapadeas, 2024), investments in environmental disasters (Niu and Zou, 2024), and the optimization of model distortion subject to moment constraints (Jaimungal et al., 2024); however, its application to Jacobi processes, particularly for modeling sustainable tourism, has not yet been found. We show that even when MFG is subject to model misspecification, it admits a closed-form solution; namely, both the HJB and FP equations are solved explicitly in conjunction with a consistency condition representing the mean field effect. The unique solvability of the consistency equation follows from a fixed-point argument. The resulting controlled process is a mean field Jacobi process, which is a Jacobi-type process whose behavior is classified according to its parameter values so that sustainable tourism and overtourism states can be characterized analytically without simulating sample paths. This tractability of the controlled process is advantageous because numerically simulating an SDE with degenerate diffusion is difficult (Hefter and Herzwurm, 2018).

The third contribution of this study is the detailed investigation of the MFG for cases in which a closed-form solution is unavailable. We propose a finite difference method to overcome this difficulty, where the stationary HJB and FP equations are solved through a forward–forward coupling (Achdou and Capuzzo-Dolcetta, 2010; Festa et al., 2024). The combined use of a monotone numerical Hamiltonian in the discretized HJB equation and an upwind method in the FP equation enables the acquisition of physically consistent numerical solutions. Similar numerical methods were successfully applied to other MFG models (Ráfales and Vázquez, 2021; Ráfales and Vázquez, 2024). The finite difference method is verified against a closed-form solution, after which more complex cases are investigated using this method.

### 1.3 Organization of this paper

The remainder of this paper is organized as follows. **Section 2** reviews the classical Jacobi process and presents its mean field version. The notion of model uncertainty applied to the Jacobi process is introduced in **Section 2**. **Section 3** formulates an MFG for sustainable tourism and derives a closed-form solution under certain conditions. **Section 4** is devoted to the application of the MFG. **Section 5** concludes the paper and presents future perspectives. **Appendices** contain auxiliary results, such as proofs of propositions and the discretization scheme of the finite difference method.

## 2. Jacobi processes

Throughout this study, we use a complete probability space $(\Omega, \mathbb{F}, \mathbb{P})$ ($\Omega$ is collection of all events, $\mathbb{F}$ is filtration, and $\mathbb{P}$ is probability function), as in conventional stochastic calculus (e.g., Definition 2.1. in Capasso and Bakstein (2021)).

### 2.1 Classical Jacobi process

We set the domain $D = (0,1)$ and its closure $\bar{D} = [0,1]$. The classical Jacobi process $X = (X_t)_{t \geq 0}$ is governed by the Itô's SDE (Chapter 6 in Alfonsi (2015))

$$dX_t = \underbrace{(a - bX_t)dt}_{\text{Mean reversion}} + \underbrace{\sigma\sqrt{X_t(1 - X_t)}dB_t}_{\text{Stochastic fluctuation}}, \ t > 0 \quad (1)$$

subject to an initial condition $X_0 \in D$, where $a > 0$ is the source parameter, $b > 0$ is the reversion parameter, $\sigma > 0$ is the volatility, and $B = (B_t)_{t \geq 0}$ is a 1-D standard Brownian motion. We set the stopping times $\tau_0$ and $\tau_1$:

$$\tau_0 = \inf\{t > 0 | X_t = 0\} \text{ and } \tau_1 = \inf\{t > 0 | X_t = 1\}. \tag{2}$$

***Remark 1*** In the context of the tourism modeling in this paper, $X$ corresponds to the tourism demand, which is a conceptual variable such that an agent visits a target tourist site at a higher frequency with a larger value of $X$. The solution $X$ can therefore be interpreted as a normalized arrival frequency or its increasing function of an agent at the tourist site. The state $X_t = 0$ implies that the agent does not visit the site.

According to Theorem 6.1.1 of Alfonsi (2015), the SDE (1) admits a unique pathwise solution that is continuous and bounded in $\bar{D}$ almost surely (a.s.) $t > 0$ under **Assumption 1**.

***Assumption 1***

$$a \leq b. \tag{3}$$

Moreover, according to Proposition 6.1.2 of Alfonsi (2015), the SDE (1) admits a unique pathwise solution that is continuous and bounded in $D$ a.s. for $t > 0$ ($\tau_0 = \tau_1 = +\infty$) under **Assumption 2**.

***Assumption 2***

$$2a \geq \sigma^2 \text{ and } 2(b-a) \geq \sigma^2 > 0. \tag{4}$$

Here, **Assumption 1** implies that the mean reversion is sufficiently strong. **Assumption 2** implies that the mean reversion is sufficiently strong, and further the noise intensity is sufficiently small. **Assumption 1** is satisfied under **Assumption 2**.

The FP equation that governs the stationary PDF $p = p(x)$ of the Jacobi process, namely the SDE (1), is given as follows (Example 5.46 in Capasso and Bakstein (2021))

$$\frac{\mathrm{d}}{\mathrm{d}x}\left((a-bx)p - \frac{\mathrm{d}}{\mathrm{d}x}\left(\frac{\sigma^2}{2}x(1-x)p\right)\right) = 0, \quad x \in D \tag{5}$$

along with the boundary condition

$$(a-bx)p - \frac{\mathrm{d}}{\mathrm{d}x}\left(\frac{\sigma^2}{2}x(1-x)p\right) = 0, \quad x = 0,1. \tag{6}$$

Under **Assumption 2**, the FP equation is explicitly solved as

$$p(x) = \frac{x^{\alpha-1}(1-x)^{\beta-1}}{\mathrm{Be}(\alpha,\beta)}, \quad x \in D, \tag{7}$$

where $\alpha = \frac{2a}{\sigma^2} > 0$ and $\beta = \frac{2(b-a)}{\sigma^2} > 0$. Here, $\mathrm{Be}(\alpha,\beta)$ is the beta function given by

$$\text{Be}(\alpha,\beta) = \int_0^1 x^{\alpha-1}(1-x)^{\beta-1}\,\mathrm{d}x. \tag{8}$$

The stationary average of the Jacobi process is $\mathbb{E}[X_t] = \dfrac{a}{b} \in \bar{D}$.

If $\beta = 0$, then $a = b$ and we have $\tau_1 < +\infty$ and $X_t = 1$ for $t > \tau_1$, showing that the stationary distribution is Dirac delta concentrated at $x = 1$ and $p(x) = \delta(x-1)$. Similarly, for $a = 0$ and $b > 0$, the stationary distribution is the Dirac delta concentrated at $x = 0$ ($p(x) = \delta(x)$). If $a > b$, then $\tau_1 < +\infty$; additionally, we set $X_t = 1$ if $t > \tau_1$; otherwise, the SDE cannot be defined for $t > \tau_1$. This is because the drift coefficient becomes positive at $X_{\tau_1} = 1$ ($a - bX_{\tau_1} > 0$) and is hence upward, whereas the diffusion coefficient $\sqrt{X_t(1-X_t)}$ is defined only if $X_t \in [0,1]$ and becomes undefinable after $\tau_1$.

*Remark 2* Considering **Remark 1**, in the context of modeling sustainable tourism, the state $X_t = 0$ for all $t \geq 0$ corresponds to no tourist at the target site, whereas the state $X_t = 1$ for all $t \geq 0$ corresponds to the overtourism state at which the site would be significantly damaged. The state $X_t \in D$ for all $t \geq 0$ corresponds to internal dynamics, suggesting the establishment of sustainable tourism at the target site. The same classification applies to the mean field version introduced in **Section 2.2**. Under **Assumption 2**, it follows that $p(0) = p(1) = 0$ because $\alpha > 1$ and $\beta > 1$ (sustainable tourism state). The PDF satisfies $p(1) > 0$ if $\beta = 1$ (again a sustainable tourism state) and is unbounded as $p(1) = +\infty$ if $0 < \beta < 1$ (coexistence between sustainable tourism and overtourism states). The stationary PDF (7) is not definable if $\beta \leq 0$ i.e., $b \leq a$ (overtourism state).

**2.2 Mean field Jacobi process**

A mean field SDE, also referred to as a McKean–Vlasov SDE, is an Itô's SDE whose drift or diffusion coefficient depends on the law, e.g., PDF, of the process itself. Our mean field Jacobi process has the form

$$\mathrm{d}X_t = \left(\hat{a}(\bar{X}_t) - \hat{b}(\bar{X}_t)X_t\right)\mathrm{d}t + \sigma\sqrt{X_t(1-X_t)}\,\mathrm{d}B_t, \; t > 0 \tag{9}$$

subject to an initial condition $X_0 \in D$, where $\hat{a}: \mathbb{R} \to [0,+\infty)$ and $\hat{b}: \mathbb{R} \to [0,+\infty)$ are source and reversion parameters that depend on the mean $\bar{X}_t = \mathbb{E}[X_t]$ of the process $X$. The process $X$ here can be perceived as the tourism demand of a representative agent among many other ones following the same form of SDEs.

To the best of the author's knowledge, although the well-posedness of the process is theoretically covered by the general theory of McKean–Vlasov SDEs, as discussed below, the mean field version of the Jacobi process has not been explicitly dealt with in the literature. In this subsection, we assume **Assumptions 3** and **4** to well-pose the mean field SDE (9).

***Assumption 3*** *The coefficients $\hat{a}$ and $\hat{b}$ are strictly bounded and satisfy the following global Lipschitz continuity condition: For some constant $K > 0$, it follows that*

$$|\hat{a}(x_1) - \hat{a}(x_2)|, |\hat{b}(x_1) - \hat{b}(x_2)| \leq K|x_1 - x_2|, \quad x_1, x_2 \in \mathbb{R}. \tag{10}$$

***Assumption 4*** *The coefficients $\hat{a}$ and $\hat{b}$ satisfy*

$$\hat{a}(x_1) \leq \hat{b}(x_1), \quad x_1 \in \mathbb{R}. \tag{11}$$

Through direct application of Theorem 3.1. of Liu et al. (2023) under **Assumptions 3** and **4**, the modified SDE, which is SDE (9) with the diffusion coefficient replaced by $\sigma\sqrt{|X_t(1-X_t)|}$, admits a unique pathwise solution that is continuous in time. Moreover, by assuming an establishment of a stationary state at time 0, applying the contradiction argument (e.g., the argument used in Theorem 6.1.1 of Alfonsi (2015)) to this modified SDE shows that its unique solution coincides with that of the mean field SDE (9).

We provide a multiagent interpretation of the mean field SDE (9) so that its background can be better understood. The explanation below is based on Cardaliaguet and Porretta (2021). The mean field SDE (9) can be considered as a large population limit of the following $N$-dimensional multiagent system ($N \in \mathbb{N}$) that governs the $N$-dimensional process $(X_{1,t}, X_{2,t}, ..., X_{N,t})_{t \geq 0}$:

$$dX_{i,t} = \left(\hat{a}\left(\frac{1}{N-1}\sum_{j \neq i}^{N} X_{j,t}\right) - \hat{b}\left(\frac{1}{N-1}\sum_{j \neq i}^{N} X_{j,t}\right)X_{i,t}\right)dt + \sigma\sqrt{X_{i,t}(1-X_{i,t})}dB_{i,t}, \quad t > 0, \tag{12}$$

where $(B_{1,t}, B_{2,t}, ..., B_{N,t})_{t \geq 0}$ is a collection of $N$-dimensional independent standard one-dimensional Brownian motions. Taking the limit $N \to +\infty$ in the system (12) yields the convergence $\frac{1}{N-1}\sum_{j \neq i}^{N} X_{j,t} \to \bar{X}_t$ owing to the law of large numbers, and all $X_i$ follow the SDE of the form (9).

Our formulation directly deals with the mean field SDE because our interest is not in the convergence from finite to infinite players but rather in the effective characterization of macroscopic tourism states. We note that $x$-dependent version of **Assumption 2** is actually unnecessary for this manuscript because our model turns out to be a Jacobi diffusion process where $a, b$ in **Assumption 2** are constant as shown in **Section 3.3**. Note that by **Assumption 4** solutions to (9) will be in $\bar{D}$ but not necessarily in $D$.

## 2.3 Model uncertainty

Misspecification, namely model uncertainty, is described as a distortion of the Brownian motion $B$ in the sense of Knight (Hansen and Sargent, 2001). Specifically, the probability measure on which the SDE is defined is denoted as $\mathbb{P}$. We then consider a real-valued measurable process $w = (w_t)_{t \geq 0}$ such that there exists a probability measure $\mathbb{Q} = \mathbb{Q}(w)$ (depending on $w$) such that the following process $\hat{B} = (\hat{B}_t)_{t \geq 0}$ becomes a standard Brownian motion under the alternative probability measure $\mathbb{Q}$:

$$\hat{B}_t = B_t - \int_0^t w_s \mathrm{d}s, \ t > 0, \tag{13}$$

with which the mean field SDE (9) under $\mathbb{Q}$ becomes

$$\mathrm{d}X_t = \left(\hat{a}(\bar{X}_t) - \hat{b}(\bar{X}_t)X_t + \underbrace{\sigma\sqrt{X_t(1-X_t)}w_t}_{\text{Additional drift}}\right)\mathrm{d}t + \sigma\sqrt{X_t(1-X_t)}\mathrm{d}\hat{B}_t, \ t > 0, \tag{14}$$

which is a mean field Jacobi process with additional drift due to distortion. We formally express $\mathbb{P}$ according to $\mathbb{P} = \mathbb{Q}(0)$, because $\mathbb{P}$ coincides with $\mathbb{Q}(w)$ when there is no misspecification ($w_t = 0$, $t \geq 0$). The variable $w$ denotes uncertainty, and its controller is nature, a virtual decision-maker which controls the uncertainty such that it adversely affects the objective function of the agent as a traveler.

The difference between $\mathbb{P} = \mathbb{Q}(0)$ and $\mathbb{Q} = \mathbb{Q}(w)$ is evaluated using the relative entropy $\mathbb{D}(\mathbb{Q}|\mathbb{P})$. For each time interval $[0,T]$ with $T > 0$, the relative entropy is expressed as (Hansen and Sargent, 2001)

$$\mathbb{D}(\mathbb{Q}|\mathbb{P}) = \mathbb{E}_{\mathbb{Q}}\left[\int_0^T \frac{w_s^2}{2}\mathrm{d}s\right]. \tag{15}$$

As implied in (15), $\mathbb{D}(\mathbb{Q}|\mathbb{P}) \geq 0$ and $\mathbb{D}(\mathbb{Q}|\mathbb{P}) = 0$ if and only if $w \equiv 0$, meaning that only the positivity of the relative entropy implies the existence of misspecification. Moreover, (15) suggests that the unit-time relative entropy at time $t$ is $\frac{w_t^2}{2}$.

## 3. Mean field game
### 3.1 Controlled system and objective function
#### 3.1.1 Controlled system

We consider an MFG to be a zero-sum stochastic differential game between a representative agent and nature, representing misspecification. We assume that a representative agent has tourism demand at a target site $X$, which is governed by the following controlled SDE on $\mathbb{Q}(w)$:

$$\mathrm{d}X_t = \left(u_t - bX_t + \sigma\sqrt{X_t(1-X_t)}w_t\right)\mathrm{d}t + \sigma\sqrt{X_t(1-X_t)}\mathrm{d}\hat{B}_t, \ t > 0 \tag{16}$$

subject to an initial condition $X_0 \in D$ ($b, \sigma > 0$ are constants), where the real-valued measurable process $u = (u_t)_{t \geq 0}$ represents the drift to modulate the tourism demand, and the real-valued measurable process $w = (w_t)_{t \geq 0}$ represents the distortion, or equivalently, the misspecification. The processes $u$ and $w$ are assumed to be controlled by the representative agent and nature, respectively. At this stage, the SDE (16) involves no mean field interaction. The admissible sets $\mathbb{U}$ and $\mathbb{W}$ of $u$ and $w$ are respectively defined as follows (measurability means with respect to a natural filtration generated by $\hat{B}$):

$$\mathbb{U} = \left\{ u = (u_t)_{t \geq 0} \middle| u_t \text{ is progressively measurable w.r.t } \mathbb{F}, \delta \mathbb{E}\left[ \int_0^{+\infty} e^{-\delta s} u_s^2 \mathrm{d}s \right] < +\infty \text{ for each } \delta > 0 \right\} \quad (17)$$

and

$$\mathbb{W} = \left\{ w = (w_t)_{t \geq 0} \middle| w_t \text{ is progressively measurable w.r.t } \mathbb{F}, \delta \mathbb{E}\left[ \int_0^{+\infty} e^{-\delta s} g(X_s) w_s^2 \mathrm{d}s \right] < +\infty \text{ for each } \delta > 0 \right\}, \quad (18)$$

where $g : \bar{D} \to [0, +\infty)$ is a bounded and measurable function specified in **Section 3.1.2**. We also require that the controlled process (16) admits a unique pathwise solution that is continuous and bounded in $\bar{D}$ for $t > 0$. This assumption excludes the controlled paths that only exist locally over time (see the discussion immediately after **Proposition 1** in **Section 3.3.1**). Controls should be non-anticipating as commonly assumed in differential games (e.g., Chapter VIII in Bardi and Capuzzo-Dolcetta (1997)), which is implicitly assumed in this paper as well.

### 3.1.2 Objective function

An objective function is the index to be maximized by the representative agent while minimized by nature. We consider an MFG in an infinite horizon, whose objective function is the mapping $J : \bar{D} \times \mathbb{U} \times \mathbb{W} \to \mathbb{R}$:

$$J(x, u, w) = \mathbb{E}_{\mathbb{Q}} \left[ \int_0^{+\infty} e^{-\delta s} \left( \underbrace{-\frac{c(\bar{X}_s)}{2} u_s^2}_{\text{Modulating cost}} + \underbrace{f(X_s)}_{\text{Tourism utility}} \right) \mathrm{d}s + \underbrace{\int_0^{+\infty} e^{-\delta s} \frac{g(X_s)}{2} w_s^2 \mathrm{d}s}_{\text{Weighted model uncertainty}} \middle| X_0 = x \right], \quad (19)$$

and its worst-case optimized version, the value function, given by

$$\Phi(x) = \sup_{u \in \mathbb{U}} \inf_{w \in \mathbb{W}} J(x, u, w), \quad x \in \bar{D}, \quad (20)$$

where $\delta > 0$ denotes the discount rate. Each term on the right-hand side of (19) can be explained as follows. The first term represents the cost of modulating travel demand, such as information collection and travel fees, with $c :\to [0, +\infty)$ as a function of the average $\bar{X}_t = \mathbb{E}_{\mathbb{Q}}[X_t]$. Later, we will consider $c$ proportional to $\bar{X}_t$ with which our mean field game becomes analytically solvable. This term conceptually represents the costs in an aggregated manner. The quadratic dependence of this term on the control $u$ is used to obtain an analytically tractable model. The weighting factor $c$ is assumed to be a continuous and nondecreasing function in $\bar{D}$ such that the cost increases when $X_t$ is larger on average. This mechanism is expected to suppress an excessively large $X_t$ potentially contributing to overtourism. The weighting factor can be designed by the local government of the target travel site by setting a tourism tax, namely, dynamic pricing (Chenavaz et al., 2022), such that it can be set higher when a large number of travelers are expected to arrive at the site. Moreover, increasing tourism demand may trigger conflicts between residents and tourists (Takahashi, 2024), which is another form of overtourism.

The second term in (19) in conjunction with the function $f : \bar{D} \to \mathbb{R}$, represents the travel utility gained by the representative agent, which is assumed to be a unimodal function maximized at some point in $D$. Travel utilities can also be designed by the local government of a target travel site to suppress overtourism. Unless otherwise specified, we assume

$$f(x) = -\frac{P}{2}x^2 + Qx, \ x \in \bar{D} \tag{21}$$

with the real parameters $P$ and $Q$ such that $P > Q$ so that $f$ has a global maximum at $x = \frac{Q}{P} \in D$.

The third term in (19) represents the weighted relative entropy between the benchmark and distorted models. The weighting factor $g$ represents the uncertainty aversion of the representative agent such that a larger value of $g$ implies that he/she assumes a smaller misspecification in the model. The third term can also be considered as a control cost to be paid by the nature; a larger $g$ implies a larger cost with which the dynamics of the SDE (16) are less affected. A nonconstant $g$ represents a situation in which the misspecification is penalized depending on the system dynamics such that a smaller value of $g = g(X_s)$ corresponds to a more pessimistic model distortion at the state $X_s$. This type of nonconstant uncertainty aversion has been employed in the literature to represent heterogeneous uncertainty aversion (Baltas, 2024; Lv et al., 2023; Mu et al., 2023).

We use the following $g$ so that the no- and overtourism states are more strongly penalized than the internal ones because these boundary states should be avoided:

$$g(x) = \eta x(1-x), \ x \in \bar{D}, \tag{22}$$

where $\eta > 0$ is a constant representing the degree of uncertainty aversion; a larger $\eta$ implies a smaller uncertainty aversion that assumes smaller misspecification, and vice versa. Another characteristic of this form of $g$ is that the optimally controlled dynamic remains a Jacobi-type SDE that is highly tractable because the associated HJB equation has the linear-quadratic form and the term singularity $\sqrt{X_t(1-X_t)}w_t$ becomes linear in $X_t$ drift under the optimality as shown in the next subsections.

### 3.2 Optimality system

Based on the dynamic programming principle assuming the establishment of a stationary state at time 0 (e.g., Calvia et al., 2024; Carmona et al., 2023; La Torre and Maggistro, 2024), the optimality system associated with our MFG contains the HJB equation

$$\begin{aligned} -\delta\Phi(x) + \sup_{u \in \mathbb{R}} \inf_{w \in \mathbb{R}} \left\{ \left(u + \sigma\sqrt{x(1-x)}w\right)\frac{d\Phi}{dx} - \frac{c(m)}{2}u^2 + \frac{\eta}{2}x(1-x)w^2 \right\} \\ -bx\frac{d\Phi}{dx} + \frac{\sigma^2}{2}x(1-x)\frac{d^2\Phi}{dx^2} + f(x) = 0 \end{aligned}, \ x \in \bar{D} \tag{23}$$

and the FP equation

$$\frac{d}{dx}\left(\left(u^*(x) - bx + \sigma w^*(x)\sqrt{x(1-x)}\right)p\right) - \frac{d}{dx}\left(\frac{\sigma^2}{2}x(1-x)p\right) = 0, \ x \in D \tag{24}$$

subject to the boundary condition

$$\left(u^*(x) - bx + \sigma w^*(x)\sqrt{x(1-x)}\right)p - \frac{d}{dx}\left(\frac{\sigma^2}{2}x(1-x)p\right), \quad x = 0,1. \tag{25}$$

Here, $m = \int_0^1 xp(x)dx$ is the stationary mean of tourism demand, and $u^*$ and $w^*$ as candidates of optimal controls are expressed through $\Phi$:

$$u^*(x) = \arg\max_{u \in \mathbb{R}} \left\{ u\frac{d\Phi}{dx} - \frac{c(m)}{2}u^2 \right\} = \frac{1}{c(m)}\frac{d\Phi}{dx} \tag{26}$$

and

$$w^*(x) = \arg\min_{w \in \mathbb{R}} \left\{ \sigma\sqrt{x(1-x)}w\frac{d\Phi}{dx} + \frac{\eta}{2}x(1-x)w^2 \right\} = -\frac{\sigma}{\eta}\frac{1}{\sqrt{x(1-x)}}\frac{d\Phi}{dx}. \tag{27}$$

The HJB equation (23) satisfies the so-called Isaacs condition that the maximization and minimization problems in the equation are decoupled. By (26) and (27), the HJB equation (23) is rewritten as

$$-\delta\Phi(x) + M(m)\left(\frac{d\Phi}{dx}\right)^2 - bx\frac{d\Phi}{dx} + \frac{\sigma^2}{2}x(1-x)\frac{d^2\Phi}{dx^2} + f(x) = 0, \quad x \in \bar{D} \tag{28}$$

and the FP equation (24) as

$$\frac{d}{dx}\left(\left(2M(m)\frac{d\Phi}{dx} - bx\right)p - \frac{d}{dx}\left(\frac{\sigma^2}{2}x(1-x)p\right)\right) = 0, \quad x \in D. \tag{29}$$

Here, the coefficient $M(m)$ is given by

$$M(m) = \frac{1}{2}\left(\frac{1}{c(m)} - \frac{\sigma^2}{\eta}\right). \tag{30}$$

The boundary condition (25) is rewritten accordingly.

The quantities $u^*$ and $w^*$ correspond to the optimal controls by the agent and nature, respectively, assuming that they are Markovian as in conventional stochastic control and mean field game models. The controlled SDE (16) then becomes

$$dX_t = \left(2M(m)\frac{d\Phi}{dx}(X_t) - bX_t\right)dt + \sigma\sqrt{X_t(1-X_t)}d\hat{B}_t, \quad t > 0 \tag{31}$$

with its stationary PDF governed by (29), suggesting that finding the value function $\Phi$ is essential for determining the controlled dynamics.

In the sequel, we always assume **Assumption 5** to verify a guessed solution for the optimal system.

***Assumption 5***

$$M(0) = \frac{1}{2}\left(\frac{1}{c(0)} - \frac{\sigma^2}{\eta}\right) > 0. \tag{32}$$

**Assumption 5** implies that there is some $m \in (0,1)$ such that $M(m) > 0$ because $c(m)$ is nondecreasing with respect to $m \in (0,1)$. Here, we set $M(0) = +\infty$ if $c(0) = 0$.

### 3.3 Closed-form solution

In this section, we consider the mean field game formulated in the previous section. We first present **Propositions 1-3** as the main outcome of this paper. Variables and parameters appearing in these propositions will be explained in the subsections below along with their backgrounds and explanations.

**Proposition 1** *Assume that there exists a unique solution $m = m^* \in D$ to the following consistency equation such that $M(m^*) > 0$:*

$$m = \frac{2M(m)E_-}{b - 2M(m)A_-}. \tag{33}$$

*Then, the value function is given by*

$$\Phi(x) = \frac{A^*}{2}x^2 + E^*x + F^*, \quad x \in \bar{D} \tag{34}$$

*and the controlled SDE at a stationary state by*

$$dX_t = \left(2M(m^*)E^* - \left(b - 2M(m^*)A^*\right)X_t\right)dt + \sigma\sqrt{X_t(1-X_t)}d\hat{B}_t, \quad t > 0. \tag{35}$$

*Here, we set $A^* = A_-\big|_{m=m^*}$, $E^* = E_-\big|_{m=m^*}$, and*

$$E_- = \frac{Q + \frac{\sigma^2}{2}A_-}{\delta + b - 2M(m)A_-} \quad \text{and} \quad F^* = \frac{M(m^*)(E^*)^2}{\delta}. \tag{36}$$

*The controlled process $X$ at stationary state follows a Jacobi process that is valued in $\bar{D}$. Finally, the optimal controls $u^*$ and $w^*$ are given by*

$$u^*(x) = \frac{1}{c(m^*)}(A^*x + E^*) \quad \text{and} \quad w^*(x) = -\frac{\sigma}{\eta}\frac{1}{\sqrt{x(1-x)}}(A^*x + E^*). \tag{37}$$

**Proposition 2** *The consistency equation (33) admits a unique solution in $D$ if*

$$\left.\frac{2M(m)E_-}{b - 2M(m)A_-}\right|_{m=0} > 0 \quad \text{and} \quad \left.\frac{2M(m)E_-}{b - 2M(m)A_-}\right|_{m=1} < 1 \tag{38}$$

*and the right-hand side of (33) is strictly decreasing with respect to $m \in \bar{D}$.*

**Proposition 3** *Assume that $M(1) > 0$ and $c(m)$ is strictly increasing with respect to $m \in \bar{D}$. In addition, assume that $P \geq 0$ is sufficiently small and the condition (56) is satisfied. Then, the consistency equation (33) admits a unique solution $m = m^*$ in $D$.*

### 3.3.1 Derivation of the solution

We derive a closed-form solution for the optimality system (28)-(29), and verify that the solution yields the value function of (20). We guess the following solution to the HJB equation (28):

$$\Phi(x) = \frac{A}{2}x^2 + Ex + F, \quad x \in \bar{D} \tag{39}$$

with constants $A, E, F \in \mathbb{R}$ to be determined. Substituting (39) into (28) yields

$$\begin{aligned}&-\delta\left(\frac{A}{2}x^2 + Ex + F\right) + M(m)\left(A^2 x^2 + 2AEx + E^2\right)\\&-bx(Ax+E) + \frac{\sigma^2}{2}(x - x^2)A - \frac{P}{2}x^2 + Qx = 0\end{aligned}, \quad x \in \bar{D}. \tag{40}$$

Given one $m \in D$, we obtain the following equations to determine $A, E, F$:

$$M(m)A^2 - \left(b + \frac{\sigma^2}{2} + \frac{\delta}{2}\right)A - \frac{P}{2} = 0, \tag{41}$$

$$(\delta + b - 2M(m)A)E = Q + \frac{\sigma^2}{2}A, \tag{42}$$

$$F = \frac{M(m)E^2}{\delta}. \tag{43}$$

Once a suitable $A$ is obtained from (41), the other constants $E, F$ can be found by using (42) and (43) provided $\delta + b - 2M(m)A \neq 0$.

We would have two solutions $A_\pm$ to the quadratic equation (41). If $M(m) > 0$, then

$$A_\pm = \frac{1}{2M(m)}\left(b + \frac{\sigma^2}{2} + \frac{\delta}{2} \pm \sqrt{\left(b + \frac{\sigma^2}{2} + \frac{\delta}{2}\right)^2 + 2M(m)P}\right). \tag{44}$$

They are real and satisfy $A_- < 0 < A_+$. We present a technical lemma, which directly follows from (44).

**Lemma 1** *If $M(m) > 0$, then*

$$b - 2M(m)A_- > 0. \tag{45}$$

**Proposition 1** is then proven based on **Lemma 1**.

An important contribution of **Proposition 1** is that the controlled process follows a Jacobi process, which is therefore analytically tractable. The sufficiently high smoothness of the value function (34) enables its verification without resorting to the notion of a weak solution, such as a viscous solution. The nonlinear equation (33) is the consistency equation to be satisfied by the average of the controlled process (e.g., see Section 2.3 in Achdou et al., 2023; Section 5.2 in Gomes et al., 2023; Proof of Lemma 3.2 in Liang and Zhang, 2024). An intuitive explanation for using $A_-$ rather than $A_+$ is that we have

$$\delta + b - 2M(m)A_{+} < 0, \tag{46}$$

and further the net reversion rate $b - 2M(m)A_{+}$ of the controlled process $X$ becomes negative and cannot be defined globally in time. Control based on $A_{+}$ is not admissible owing to the local-in-time existence of the controlled process. As discussed above, the system (28)-(29) has at least two solutions. However, the solution of the quadratic form with an admissible control is unique; we have already written in **Section 3.1.3** that the admissibility requires the global existence of a unique bounded solution to the controlled system. So, the non-global existence with $A = A_{+}$ suffices to exclude it. Although this does not state that **Proposition 1** gives the unique solution to the mean filed game, it is a reasonable one considering the linear-quadratic form of the problem. Moreover, the solution with $A = A_{-}$ is computationally obtained by the finite difference method (see **Appendix**).

The consistency equation (33) is obtained by equating the average $m$ and its alternative representation $2M(m)E_{-}/(b - 2M(m)A_{-})$ (see **Section 2.1**). Note that both $A_{-}$ and $E_{-}$ depend on $m$. **Proposition 2** then provides a sufficient condition such that (33) admits a unique solution in $D$. The proof is omitted here because it directly follows from the classical intermediate value theorem.

The left inequality of (38) is automatically satisfied due to (45) and (36). We assume that the function $c(m)$ is bounded and strictly increasing with respect to $m \in \bar{D}$ and $M(1) > 0$. Now, the right-hand side of the consistency equation (33) becomes

$$\frac{2M(m)E_{-}}{b - 2M(m)A_{-}} = \frac{2M(m)\left(Q + \frac{\sigma^2}{2}A_{-}\right)}{(\delta + b - 2M(m)A_{-})(b - 2M(m)A_{-})}. \tag{47}$$

For simplicity, we set $h = (\sigma^2 + \delta)/2$. We have

$$b - 2M(m)A_{-} = \sqrt{(b+h)^2 + 2M(m)P} - h, \tag{48}$$

and hence

$$\frac{2M(m)E_{-}}{b - 2M(m)A_{-}} = \frac{2M(m)Q - \frac{\sigma^2}{2}\frac{2M(m)P}{\sqrt{(b+h)^2 + 2M(m)P} + b + h}}{\left(\delta + \sqrt{(b+h)^2 + 2M(m)P} - h\right)\left(\sqrt{(b+h)^2 + 2M(m)P} - h\right)}. \tag{49}$$

By (49), the consistency equation (33) is compactly rewritten as

$$m = G(M(m)), \ m \in \bar{D} \tag{50}$$

with

$$G(n) = \frac{2nQ - \frac{\sigma^2}{2}\frac{2nP}{\sqrt{(b+h)^2 + 2nP} + b + h}}{\left(\delta + \sqrt{(b+h)^2 + 2Pn} - h\right)\left(\sqrt{(b+h)^2 + 2Pn} - h\right)}, \ n > 0. \tag{51}$$

The right-hand side of (51) is continuously differentiable with respect to $P \geq 0$. We obtain

$$G_{P=0}(n) = \frac{2nQ}{(\delta+b)b}, \quad n > 0. \tag{52}$$

Therefore, the consistency equation at $P = 0$ then becomes

$$m = \frac{2Q}{(\delta+b)b} M(m), \quad m \in \bar{D} \tag{53}$$

whose right-hand side is strictly decreasing with respect to $m \in \bar{D}$. The right inequality of (38) in the present case becomes

$$\frac{2Q}{(\delta+b)b} M(1) < 1 \tag{54}$$

or equivalently

$$\underbrace{(\delta+b)b}_{\text{Discount and mean reversion}} > \underbrace{Q}_{\text{Utility}} \times \left( \underbrace{\frac{1}{c(1)}}_{\text{Cost}} - \underbrace{\frac{\sigma^2}{\eta}}_{\text{Uncertainty and noise}} \right). \tag{55}$$

This condition is satisfied if the mean reversion $b$ is sufficiently large (the fluctuation in tourism demand is suppressed in a short time interval) or the discount $\delta$ is sufficiently large (the representative agent is myopic about forecasts of tourism demand). For fixed $b$ and $\delta$, the condition expressed by (55) is satisfied if the travel utility $Q$ is sufficiently small (the tourist spot is not extremely attractive), cost $c(1)$ in the overtourism state is sufficiently large (the representative agents do not prefer crowds or pollution at the tourist stie), volatility $\sigma$ (the noise intensity) is sufficiently large (tourism demand is highly variable at the time I), or the uncertainty aversion $\eta$ is small (the model uncertainty is assumed to be not very large).

Under (55), the consistency equation (33) at $P = 0$ admits a unique solution $m = m^*$ in $D$. For a sufficiently small $P > 0$, the right-hand side of (33) (being continuous with respect to $P \geq 0$) strictly decreases with respect to $m \in \bar{D}$. In this case, the consistency equation admits a unique solution $m = m^*$ in $D$ if

$$G(M(1)) < 1. \tag{56}$$

Again, this condition is satisfied if $b$ or $\delta$ is sufficiently large. With the help of the discussion above, we obtain **Proposition 3**.

***Remark 3*** If $c$ is a positive constant, meaning that there is no mean field effect, then the consistency equation (33) admits a unique solution $m = m^*$ in $D$ if $G(M(0)) \in (0,1)$ ($M(0) = M(1)$ in this case).

### 3.3.2 Implications of the solution

We discuss more on the implications of the closed-form solution of **Proposition 1**. Below, we focus on the situation where $\sigma$ is sufficiently small so that the no-tourism state is not established as a stationary state; that is, we assume that the first inequality of (57) is satisfied.

This solution corresponds to a sustainable tourism state or an overtourism state because the controlled process of $X$ is valued in $\bar{D}$. Given **Assumption 2** and **Remark 2**, a sustainable tourism state, namely the controlled process $X$ valued in $D$, is attained if

$$\left[2M(m)E_-\right]_{m=m^*} \geq \frac{\sigma^2}{2} \text{ and } \left[b - 2M(m)A_- - 2M(m)E_-\right]_{m=m^*} \geq \frac{\sigma^2}{2}. \tag{57}$$

Considering (33), (57) can be rewritten as

$$\left[2M(m)E_-\right]_{m=m^*} \geq \frac{\sigma^2}{2} \text{ and } (1-m^*)\left[b - 2M(m)A_-\right]_{m=m^*} \geq \frac{\sigma^2}{2}. \tag{58}$$

Similarly, both the sustainable tourism and overtourism states coexist if

$$\left[2M(m)E_-\right]_{m=m^*} \geq \frac{\sigma^2}{2} \text{ and } 0 < (1-m^*)\left[b - 2M(m)A_-\right]_{m=m^*} < \frac{\sigma^2}{2}. \tag{59}$$

Only the overtourism state is established if

$$\left[2M(m)E_-\right]_{m=m^*} \geq \frac{\sigma^2}{2} \text{ and } (1-m^*)\left[b - 2M(m)A_-\right]_{m=m^*} \leq 0. \tag{60}$$

In (60), $m^*$ is a understood as the solution to the consistency equation (33) in $(0,+\infty)$.

To obtain more explicit conditions for the sustainable tourism and overtourism states, we consider the cases of $P=0$, $Q>0$, and $c(m) = c_1 m$ with $c_1 > 0$. Subsequently, the consistency equation (33) becomes

$$m = \frac{Q}{(\delta+b)b}\left(\frac{1}{c_1 m} - \frac{\sigma^2}{\eta}\right), \ m \in \bar{D}, \tag{61}$$

which can be recast as the quadratic equation

$$\frac{(\delta+b)b}{Q}m^2 + \frac{\sigma^2}{\eta}m - \frac{1}{c_1} = 0, \ m \in \bar{D}. \tag{62}$$

This quadratic equation always admits two distinct real solutions $m_-$ and $m_+$ with $m_- < 0 < m_+$:

$$m_\pm = \frac{Q}{2(\delta+b)b}\left(-\frac{\sigma^2}{\eta} \pm \sqrt{\left(\frac{\sigma^2}{\eta}\right)^2 + 4\frac{(\delta+b)b}{Qc_1}}\right), \tag{63}$$

where the candidate of $m^*$ is $m_+$. An elementary calculation shows $m^* = m_+ \in (0,1)$ under (55) with $c(1)$ is replaced by $c_1$. The condition (58) for sustainable tourism becomes

$$\frac{2M(m_+)Q}{\delta+b} \geq \frac{\sigma^2}{2} \text{ and } b(1-m_+) \geq \frac{\sigma^2}{2}, \tag{64}$$

or equivalently

$$\frac{\sigma^2}{2b} \leq m_+ \leq 1 - \frac{\sigma^2}{2b}. \tag{65}$$

The condition (65) suggests that sustainable tourism will be established when the volatility $\sigma$ as noise intensity is sufficiently small, mean reversion is sufficiently large, and $m_+$ is in a moderate range. Moreover, because

$$m_+ = \frac{2}{c_1} \frac{1}{\frac{\sigma^2}{\eta} + \sqrt{\left(\frac{\sigma^2}{\eta}\right)^2 + 4\frac{(\delta+b)b}{Qc_1}}} = \begin{cases} 0 & (\eta \to +0 : \text{Large uncertainty limit}) \\ \sqrt{\frac{Q}{c_1(\delta+b)b}} & (\eta \to +\infty : \text{No uncertainty limit}) \end{cases} \quad (66)$$

and the right-hand side of (66) is increasing with respect to $\eta > 0$, misspecification leads to a smaller average tourism demand to avoid possible overtourism.

From (59), we also obtain the condition for the coexistence of both sustainable tourism and overtourism states when

$$m_+ \geq \frac{\sigma^2}{2b} \quad \text{and} \quad m_+ > 1 - \frac{\sigma^2}{2b} \quad (67)$$

and the overtourism state is established when

$$m_+ \geq \frac{\sigma^2}{2b} \quad \text{and} \quad m_+ \geq 1. \quad (68)$$

Therefore, if $\sigma$ is sufficiently small or the mean reversion is sufficiently large, then overtourism is established when $m_+ \geq 1$, or equivalently when

$$\underbrace{(\delta+b)b}_{\text{Discount and mean reversion}} \leq \underbrace{Q}_{\text{Utility}} \times \left( \underbrace{\frac{1}{c_1}}_{\text{Cost}} - \underbrace{\frac{\sigma^2}{\eta}}_{\text{Uncertainty and noise}} \right), \quad (69)$$

which resembles (55).

### 3.3.3 Ergodic control case

Finally, we discuss the long-run limit using the vanishing discount rate $\delta \to +0$, which corresponds to the ergodic control problem. This limit corresponds to the least myopic perception limit of the representative agent and nature, where the objective function (19) is reinterpreted as the time average:

$$\begin{aligned} \delta J(x,u,w) &= \mathbb{E}_{\mathbb{Q}} \left[ \int_0^{+\infty} \delta e^{-\delta s} \left( -\frac{c(\bar{X}_s)}{2} u_s^2 + f(X_s) \right) ds + \int_0^{+\infty} \delta e^{-\delta s} \frac{g(X_s)}{2} w_s^2 ds \bigg| X_0 = x \right] \\ &\underset{\delta \to +0}{\to} \lim_{T \to +\infty} \frac{1}{T} \mathbb{E}_{\mathbb{Q}} \left[ \int_0^T \left( -\frac{c(\bar{X}_s)}{2} u_s^2 + f(X_s) \right) ds + \int_0^T \frac{g(X_s)}{2} w_s^2 ds \bigg| X_0 = x \right] \\ &= \lim_{T \to +\infty} \frac{1}{T} \mathbb{E}_{\mathbb{Q}} \left[ \int_0^T \left( -\frac{c(\bar{X}_s)}{2} u_s^2 + f(X_s) \right) ds + \int_0^T \frac{g(X_s)}{2} w_s^2 ds \right] \end{aligned} \quad (70)$$

The second line is due to the vanishing-discount approach that converts the scaled discounted integral to the time-averaged one (e.g., Bao and Tang, 2023; Cohen and Zell, 2023; Hernández-Lerma et al., 2023); additionally, the third line is from an ergodic assumption such that the information at the initial time is forgotten as time elapses, which would be satisfied if there exists a unique stationary solution to the controlled process.

To find the "maxmin" of the last line of (70) with respect to $(u,w) \in \mathbb{U} \times \mathbb{W}$ (where the integrals are understood as the time average under $\delta \to +0$), the FP equation remains the same, while the HJB equation is replaced by

$$-H + M(m)\left(\frac{d\Phi}{dx}\right)^2 - bx\frac{d\Phi}{dx} + \frac{\sigma^2}{2}x(1-x)\frac{d^2\Phi}{dx^2} + f(x) = 0, \quad x \in \bar{D} \tag{71}$$

with the optimal controls $u^*$ and $w^*$ formally given by (26) and (27), where $H \in \mathbb{R}$ is the effective Hamiltonian as the optimized objective

$$H = \sup_{u \in \mathbb{U}} \inf_{w \in \mathbb{W}} \lim_{T \to +\infty} \frac{1}{T} \mathbb{E}_{\mathbb{Q}}\left[\int_0^T \left(-\frac{c(\bar{X}_s)}{2}u_s^2 + f(X_s)\right)ds + \int_0^T \frac{g(X_s)}{2}w_s^2 ds\right]. \tag{72}$$

A solution to (71) is a combination of a constant $H \in \mathbb{R}$ and a function $\Phi: \bar{D} \to \mathbb{R}$, where $\Phi$ is an auxiliary function (rather than a value function) used to find optimal controls. Moreover, if $(H, \Phi)$ solves (71), $(H, \Phi + c')$ with an arbitrary constant $c'$ will also do so.

Now, the optimality system now contains (71), (24), and (27). Its solution is found as follows due to the guessed-solution technique, where the subscript "E" represents the ergodic control case: assume that there exists a unique solution $m = m^* \in D$ to the following equation such that $M(m^*) > 0$:

$$m = \frac{2M(m)E_E}{b - 2M(m)A_E}, \tag{73}$$

where $\Phi$ is given with an arbitrary constant $c' \in \mathbb{R}$ by

$$\Phi(x) = \frac{A_E^*}{2}x^2 + E_E^* x + c', \quad x \in \bar{D} \tag{74}$$

and the controlled SDE by

$$dX_t = \left(2M(m^*)E_E^* - (b - 2M(m^*)A_E^*)X_t\right)dt + \sigma\sqrt{X_t(1-X_t)}d\hat{B}_t, \quad t > 0. \tag{75}$$

Herein, we set $A_E^* = A_E|_{m=m^*}$, $E_E^* = E_E|_{m=m^*}$, and

$$E_E = \frac{Q + \frac{\sigma^2}{2}A_E}{b - 2M(m)A_E} \quad \text{and} \quad A_E = \frac{1}{2M(m)}\left(b + \frac{\sigma^2}{2} - \sqrt{\left(b + \frac{\sigma^2}{2}\right)^2 + 2M(m)P}\right). \tag{76}$$

Finally, the effective Hamiltonian $H$ is given by

$$H = M(m^*)(E_E^*)^2. \tag{77}$$

**Proposition 3** and **Remark 3** still apply in the ergodic control case. Hence, the classification scheme of tourism states based on a controlled SDE does.

## 4. Applications
### 4.1 Numerical method
We present a finite difference method to compute the optimality system (see (28) and (29)) for a generic

$f$, which is not necessarily quadratic. Our discretization is based on the forward–forward MFG (Achdou and Capuzzo-Dolcetta, 2010; Festa et al., 2024). It uses a monotone (upwind) discretization to compute numerical solutions stably (Ráfales and Vázquez, 2021; Ráfales and Vázquez, 2024). Choosing a sufficiently small time increment achieves stable numerical computation. A detailed discretization scheme is presented in the **Appendix**. The idea behind the proposed numerical method is that the solution to the optimality system would be obtained as a stationary solution to the following forward–forward system, subject to a suitable initial guess:

$$\frac{\partial \Phi(\tau,x)}{\partial \tau} = -\delta\Phi(\tau,x) + M(m(\tau))\left(\frac{\partial \Phi(\tau,x)}{\partial x}\right)^2 \\ -bx\frac{\partial \Phi(\tau,x)}{\partial x} + \frac{\sigma^2}{2}x(1-x)\frac{\partial^2 \Phi(\tau,x)}{\partial x^2} + f(x) \quad , \tau > 0, \ x \in \bar{D} \quad (78)$$

and

$$\frac{\partial p(\tau,x)}{\partial \tau} = -\frac{\partial}{\partial x}\left(\left(2M(m(\tau))\frac{\partial \Phi(\tau,x)}{\partial x} - bx\right)p(\tau,x) - \frac{\partial}{\partial x}\left(\frac{\sigma^2}{2}x(1-x)p(\tau,x)\right)\right), \ \tau > 0, \ x \in D \quad (79)$$

subject to the boundary condition

$$\left(2M(m(\tau))\frac{\partial \Phi(\tau,x)}{\partial x} - bx\right)p(\tau,x) - \frac{\partial}{\partial x}\left(\frac{\sigma^2}{2}x(1-x)p(\tau,x)\right) = 0, \ \tau > 0, \ x = 0,1 \quad (80)$$

along with the consistency condition

$$m(\tau) = \int_0^1 p(\tau,x)\mathrm{d}x, \ \tau > 0. \quad (81)$$

The independent variable $\tau \geq 0$ is the artificial time required to advance the system (78)–(81) such that a stationary solution is eventually obtained under the limit $\tau \to +\infty$. Accordingly, the value functions $\Phi$ and PDF $p$ are dependent on the artificial time $\tau$.

*Remark 5* Some studies have suggested that a stationary solution to a forward–forward MFG system is not necessarily attained (Gomes et al., 2016; Sun and Wang, 2024). We successfully obtained stationary solutions at least computationally.

### 4.2 Results and discussion
#### 4.2.1 Closed-form solution
We investigated the closed-form solution in **Proposition 1**. The main task in completely deriving the solution is to find the solution to the consistency equation (33). We solve this equation using the following relaxed fixed-point iteration:

$$m^{(n+1)} = (1-R)m^{(n)} + R\frac{2M(m)E_-}{b-2M(m)A_-}\bigg|^{(n)}, \ n = 0,1,2,... \quad (82)$$

subject to the initial guess that $m^{(0)} = 10^{-10}$ and a relaxation factor $R \in (0,1)$. Herein, the superscript $(n)$

represents the value computed in the $n$ th iteration. Using an excessively large $R$ results in divergence of the iteration method. In this study, we set $R = 0.01$ so that the iteration method could handle a variety of parameter values. The iteration procedure is terminated if the difference between the current and previous $m$ is smaller than $10^{-13}$.

The computed closed-form solutions are also used to validate the finite difference method, as shown in **Appendix**. In this subsection, we use parameter values and coefficients unless otherwise specified. They are chosen so that the dependence of the MFG on parameter can be visibly understood: $b = 0.1$, $\sigma = 0.2$, $P = 0.40$, $Q = 0.25$, $c(m) = c_1 m$ with $c_1 = 10$, $\eta = 1$, $\delta = 0.1$. In the present setting, the inequality (38) is satisfied automatically. Herein, we focus on $m^*$, because it determines the qualitative behavior of the controlled process.

**Figures 1-4** show the computed $m^* = m^*(P)$ for different values of the volatility $\sigma$, utility parameter $Q$, uncertainty aversion $\eta$, and discount $\delta$. Parameter $P$ is chosen as an independent parameter for visualization. This parameter serves to move the maximum point $x = \frac{Q}{P}$ of the utility $f(x)$; therefore, increasing $P$ moves the maxima to the left, representing a representative agent who has little interest in traveling. **Figures 1-4** show that a larger fluctuation of the tourism demand, smaller travel utility, stronger uncertainty aversion, or foresight with shorter perspectives reduce the average $m^*$. All the parameter dependencies visualized in these figures are monotonic.

We also investigate the regimes of sustainable tourism, the coexistence of sustainable tourism and overtourism, and overtourism. **Figure 5** shows the computed classification of tourism regimes as a function of the utility parameters $P$ and $Q$. We examine the nominal, smaller noise intensity (deterministic), larger uncertainty aversion, and smaller discount (ergodic) cases. In all cases, a sufficiently larger value of $P$ compared with $Q$ achieves a sustainable tourism state, whereas a smaller value of $P$ compared with $Q$ results in an overtourism state. This is in accordance with the quadratic shape of the utility $f$, in that its maximum is attained at $\frac{Q}{P}$, which eventually modulates the average $m^*$. The presence of noise (**Figures 5(a)-(c)** for $\sigma > 0$) leads to the appearance of a transition state in which both the overtourism and sustainable tourism states coexist; this is absent in the deterministic case (**Figure 5(a)** for $\sigma = 0$). This transition state widens under a larger uncertainty aversion of the representative agent **(Figure 5(c))**, which represents his/her hesitation about traveling owing to a possible misspecification. The ergodic case results in a wider region for achieving sustainable tourism in the $P$-$Q$ space (**Figure 5(d)**), suggesting that an agent with a longer perspective travels more actively to a tourist spot.

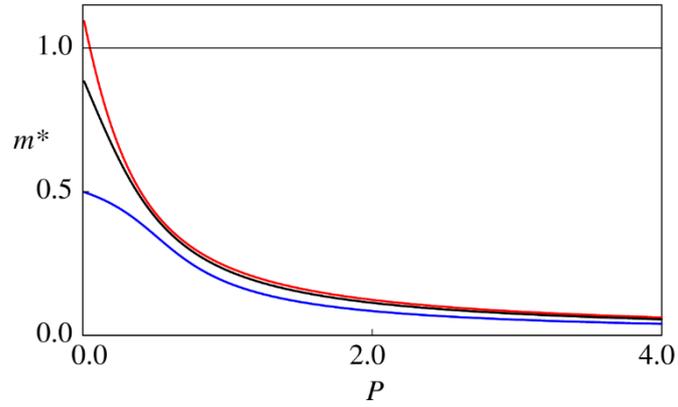

**Fig. 1** Plots of computed average $m^* = m^*(P)$ of the closed-form solutions for different values of $\sigma$ ($\sigma = 0$ (red), $\sigma = 0.2$ (black), and $\sigma = 0.4$ (blue))

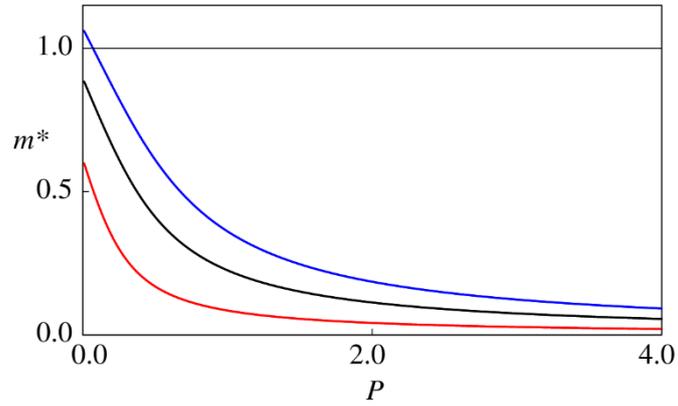

**Fig. 2** Plots of computed average $m^* = m^*(P)$ of the closed-form solutions for different values of $Q$ ($Q = 0.10$ (red), $Q = 0.25$ (black), and $Q = 0.40$ (blue))

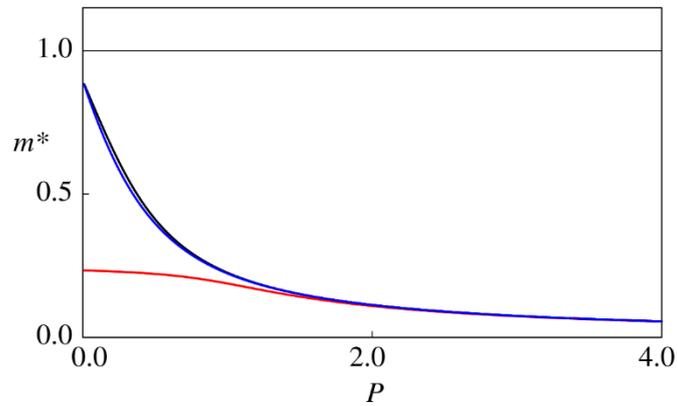

**Fig. 3** Plots of computed average $m^* = m^*(P)$ of the closed-form solutions for different values of $\eta$ ($\eta = 0.1$ (red), $\eta = 1$ (black), and $\eta = 10$ (blue))

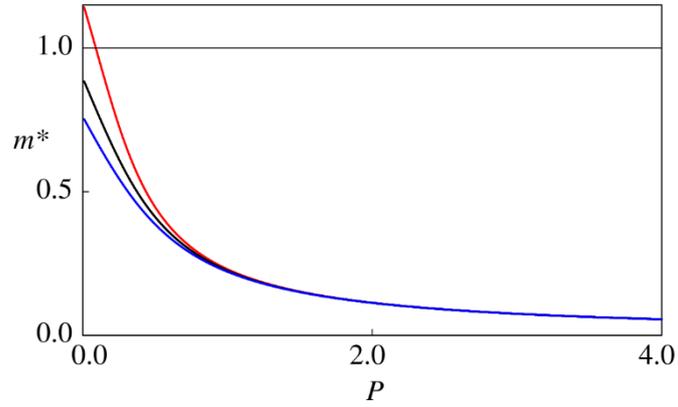

**Fig. 4** Plots of computed average $m^* = m^*(P)$ of the closed-form solutions for different values of $\delta$ ($\delta = 0$ (red), $\delta = 0.1$ (black), and $\delta = 0.2$ (blue))

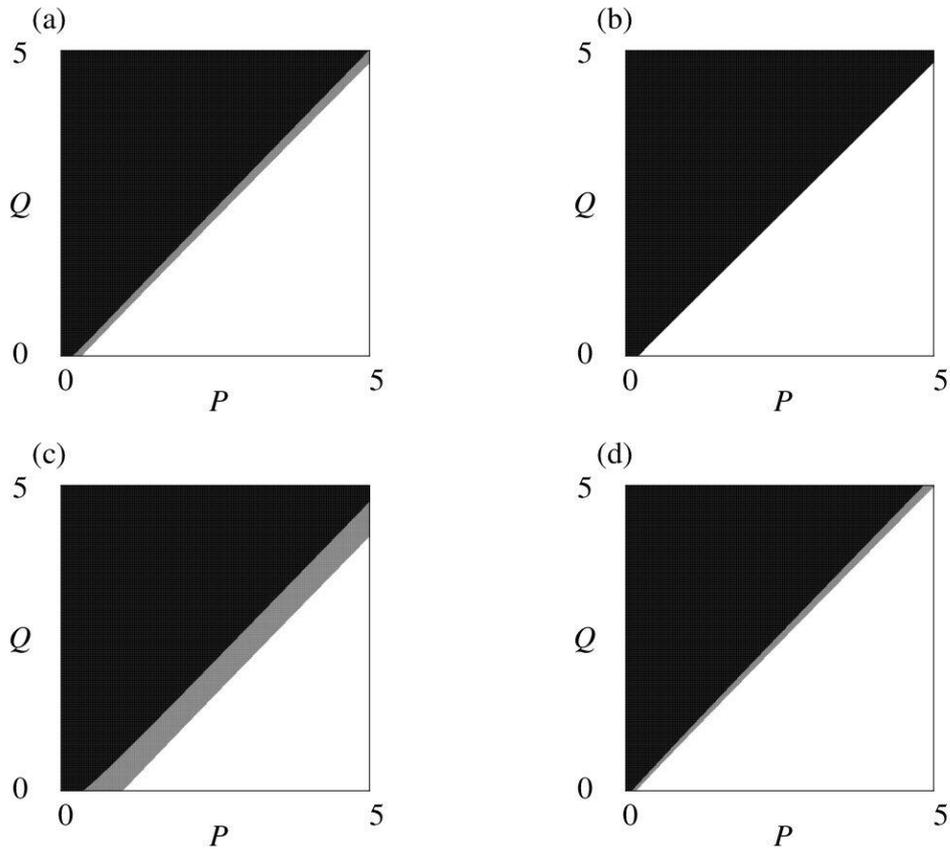

**Fig. 5** Computed classification of tourism regimes as a function of the utility parameters $P$ and $Q$, where black, gray, and white regions represent the overtourism, coexistence of sustainable tourism and overtourism, and sustainable tourism states. (a) Nominal case, (b) smaller volatility ($\sigma = 0$, deterministic case), (c) larger uncertainty aversion ($\eta = 0.5$), and (d) a smaller discount ($\delta = 0$, ergodic case)

### 4.2.2 Case with a discontinuous $f$

We apply the validated finite difference method to the solution of the MFG model for cases where closed-form solutions were not found. This case corresponds to a problem where too small or large tourist flows are not preferred, while only a moderate number of them are welcomed considering the capacity of the tourism site. We consider a piecewise constant $f$ as a simple discontinuous case. Moreover, this is the interesting case where $f$ is neither concave nor convex.

We chose the discontinuous $f$: $f(x)=1$ for $1/4 < x < 3/4$ and $f(x)=0$ otherwise, so that the travel utility is gained only for a range of moderate values of $x$. We examine both discounted ($\delta > 0$) and ergodic cases ($\delta = 0$). We use the following parameter values and coefficients unless specified otherwise; $b = 0.1$, $\sigma = 0.2$, and $c(m) = c_1 m$ with $c_1 = 5$. These have been chosen so that the dependence of the MFG on parameters can be visibly understood.

**Figures 6-7** show the computed value functions $\Phi$ for different values of the volatility $\sigma$ and discount rate $\delta$, and $\eta$. **Figures 8-9** show the computed PDFs $p$ for different values of $\sigma$ and $\delta$, and $\eta$. The computed value functions have two kinks as shown in **Figures 6-7**. These kinks, which do not appear in the closed-form solution in **Proposition 1**, are due to the discontinuity of the present $f$, demonstrating that the regularity of $f$ is inherited in $\Phi$ (e.g., see **Figures A1-A4**). **Figures 8-9** suggest that the PDFs of the controlled process are nonsmooth such that $p$ is close to zero for small $x$ values, and become positive with respect to $x$ in a nonsmooth manner owing to the kinks of the corresponding value functions in **Figures 6-7**. Increasing the uncertainty aversion with a smaller $\eta$ value shifts the PDF; hence, the average $m^*$ shifts to the left as in the closed-form solution. The PDFs are more concentrated around the center of the domain $D$ compared with the beta distributions associated with the classical Jacobi processes. These types of solutions cannot be reproduced by the closed-form solution; hence, they should be computed numerically.

**Figure 10** shows the computed average $m^*$ as a function of $\sigma$ and $\delta$ for $\eta = 1$ **(Figure 10(a))** and $\eta = 0.05$ **(Figure 10(b))**. The computed $m^*$ is respectively increasing and decreasing with respect to the examined values of $\sigma$ and $\delta$ for larger $\eta$ values, while the same observation does not hold for smaller $\eta$ values (larger uncertainty) with a small $\delta$ and a large $\sigma$ (longer perspectives with larger noise intensity). To analyze in depth this phenomenon, **Figure 11** shows the computed PDFs for $\delta = 0$ and $\eta = 0.05$, and different values of $\sigma$. **Figure 11** considers the nonmonotonic behavior of $m^*$ observed in **Figure 10(b)** and suggests that there is a sudden transition of the PDF shapes at a relatively large $\sigma$ such that $p(0)$ transits from zero to a positive value; the positive value marks the appearance of a no-tourism state, which is beyond the quadratic case where the closed-form solution is available. The transition seems to be complex according to **Figure 10(b)**; however, the finite difference method could deal with problems associated with nonquadratic $f$ with which closed-form solutions to the MFG model have not been found.

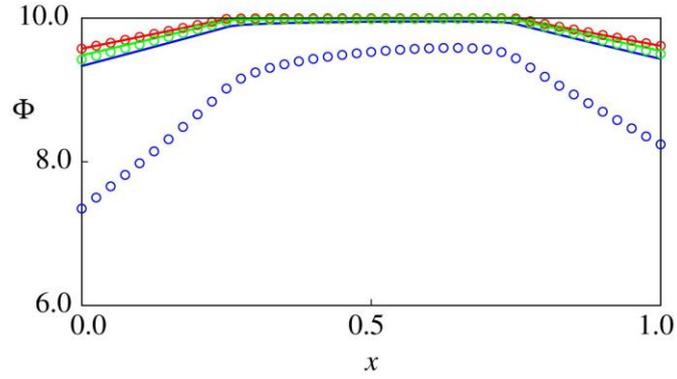

**Fig. 6** Plots of the computed value function $\Phi$ for different values of $\sigma$ and $\eta$. Symbols: $\eta = 10$ (lines) and $\eta = 0.1$ (circles). Colors: $\sigma = 0$ (red), $\sigma = 0.1$ (green), and $\sigma = 0.2$ (blue)

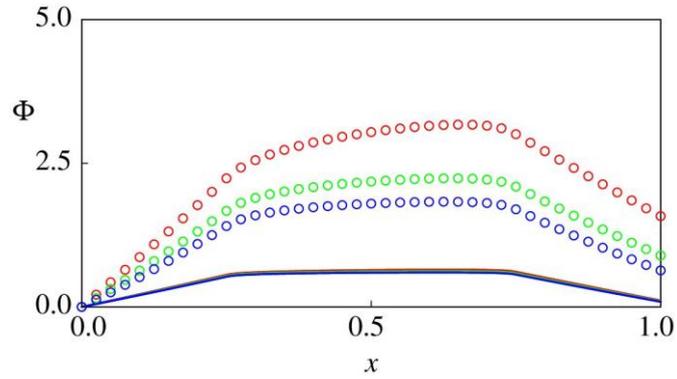

**Fig. 7** Plots of the computed value function $\Phi$ for different values of $\delta$ and $\eta$. Symbols: $\eta = 10$ (lines) and $\eta = 0.1$ (circles). Colors: $\delta = 0$ (red), $\delta = 0.1$ (green), and $\delta = 0.2$ (blue). We normalized $\Phi$ as $\Phi(0) = 0$ for visualization

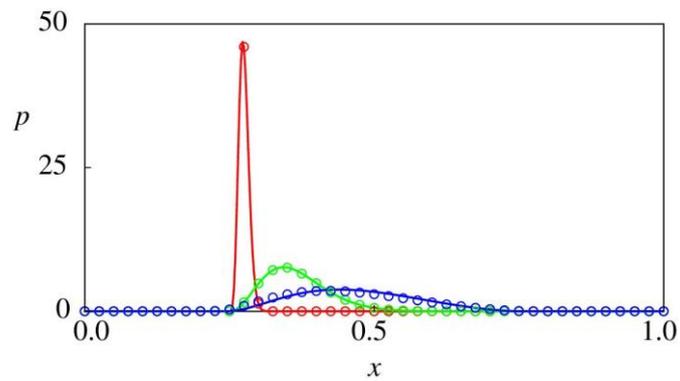

**Fig. 8** Plots of computed probability density functions (PDFs) $p$ for different values of $\sigma$ and $\eta$. Symbols: $\eta = 10$ (lines) and $\eta = 0.1$ (circles). Colors: $\sigma = 0$ (red), $\sigma = 0.1$ (green), and $\sigma = 0.2$ (blue)

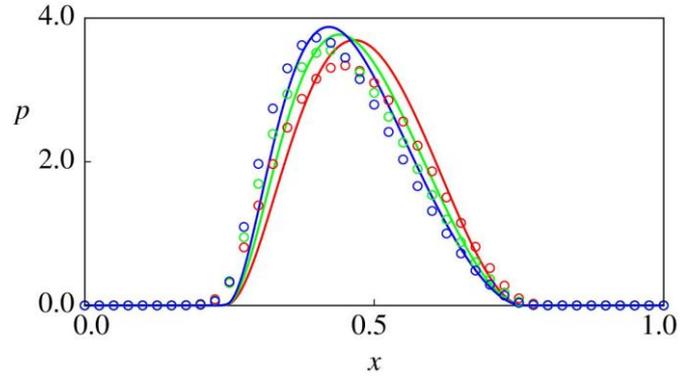

**Fig. 9** Plots of computed PDFs $p$ for different values of $\delta$ and $\eta$. Symbols: $\eta = 10$ (lines) and $\eta = 0.1$ (circles). Colors: $\delta = 0$ (red), $\delta = 0.1$ (green), and $\delta = 0.2$ (blue)

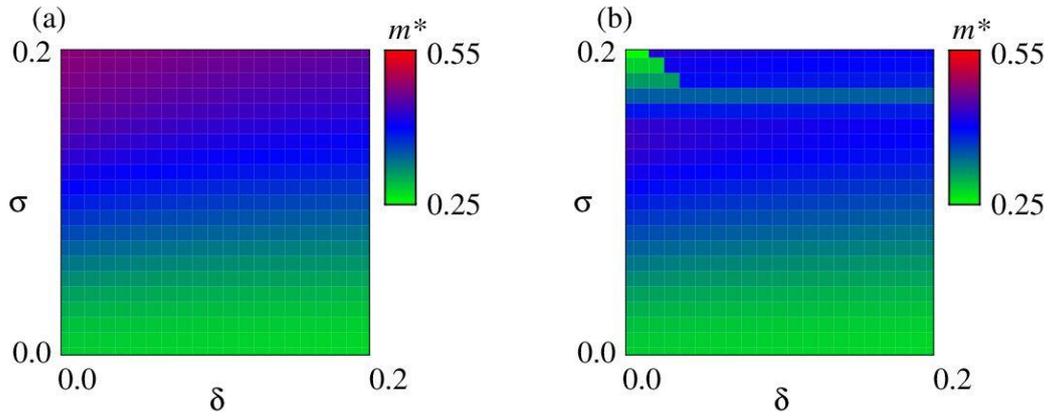

**Fig. 10** Computed average $m^*$ maps for different values of $\sigma$ and $\delta$. (a) $\eta = 1$ and (b) $\eta = 0.03$

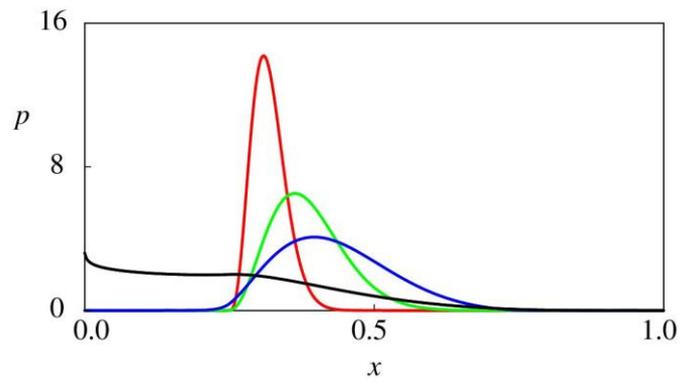

**Fig. 11** Computed PDFs for $\delta = 0$, $\eta = 0.05$, and different values of $\sigma$ ($\sigma = 0.05$ (red, $m^* = 0.319$), $\sigma = 0.10$ (green, $m^* = 0.387$), $\sigma = 0.15$ (blue, $m^* = 0.431$), and $\sigma = 0.20$ (black, $m^* = 0.250$))

**4.2.3   Implications for tourism in rural areas**

Finally, the implications of the theoretical and computational results for tourism in Japan are discussed. We discuss the implications at the end of **Section 4** because it does not contain quantitative data. Nevertheless, the model could be qualitatively applied to the target site if sufficient amount of tourism data become available, and we consider that this kind of discussion connecting mean filed games to applications is not abundant in literature.

Shiramine Village, with an area of 163.1 (km$^2$)[1] and 665 residents (as of 2024)[2], extends over the upstream reaches of the Tedori River flowing in Hakusan City, Ishikawa Prefecture in central Japan. The Tedori River is an A-class river in Japan that plays a vital role in the water resource supply in the Ishikawa Prefecture. The area surrounding the Tedori River, including Shiramine Village, was registered as a UNESCO Geopark in May 2023 because of its unique water environment, ecosystem, and geology[3]. Geoparks are expected to promote tourism and education in the future (Koizumi and Chakraborty, 2016). Shiramine Village has hot springs and has made efforts to promote tourism to sustain its local economy and society[4]. Part of the village is the Hakusan Biological Reserve, whose environment should be carefully managed. The village and its surrounding areas are famous for mountain stream fishing of the fish Iwana *Salvelinus leucomaenis* and Yamame *Oncorhynchus masou*, which are the major inland fishery resources in Japan. Inland fisheries in this area are authorized by the Hakusan-Shiramine Inland Fishery Cooperative[5].

Mammadova (2019) noted that many biological reserves, including the Hakusan Biological Reserve, face issues such as depopulation and aging, suggesting that the tourist capacity of Shiramine Village, which is the area closest to the Hakusan Biological Reserve, should be considered when designing tourism in the village. Moreover, this biological system is expected to serve as a learning site for the conservation of natural and cultural resources (Mammadova, 2017). Its role in learning is considered important because it is registered as a part of the Geopark. However, the overfishing or illegal fishing of inland fishery resources should be avoided. Indeed, illegal fishing has recently become an issue in Shiramine Village[6].

The flow of tourists to the Shiramine Village should be sustained or even increased, while the surrounding environment and ecosystems should be conserved. According to the proposed model, this implies that the sustainable tourism state should be realized in this area with some effort. For example, the village would be able to inform their characteristics globally to obtain a small uncertainty aversion corresponding to a large $\eta$. Based on the consideration listed above, overtourism in this village would result in environmental pollutions and threats to fishery resources. If the Allee effect drives the fish populations (Sass et al., 2021), for each fish species, there exists a critical population threshold; if the population falls below this threshold, the species will become extinct. The utility of the threshold type computationally examined in the previous subsection is more adequate than the parabolic one used to derive the closed-form solution. Moreover, climatic changes affect the water environment and hydrology, and habitats for the inland fishery resources, would serve as the model distortion owing to its unpredictability. Tourism policy-making at the Shiramine Village can then be addressed by incurring some environmental tax to conserve its surrounding nature. This can be incorporated into the coefficient $c(m)$ in the objective function. In view of the computational results in the previous subsections, model uncertainty should be kept sufficiently small so that the no-tourism state does not dominate. Monitoring and sharing environmental and ecological conditions of the tourist spot by a policy maker would be able to address this issue.

---

[1] https://geoshape.ex.nii.ac.jp/ka/resource/17/172103320.html (in Japanese, last accessed on August 19, 2024)
[2] https://www.city.hakusan.lg.jp/_res/projects/default_project/_page_/001/011/975/jinkou_r0601.pdf (in Japanese, last accessed on August 19, 2024)
[3] https://hakusan-geo.jp/ (in Japanese, last accessed on August 19, 2024)
[4] https://shiramine.info/index.html (in Japanese, last accessed on August 19, 2024)
[5] http://www.asagaotv.ne.jp/~gyokyou/index.html (in Japanese, last accessed on August 19, 2024)
[6] Personal communication with a union member of the Hakusan-Shiramine Inland Fishery Cooperative on August 8, 2024.

## 5. Conclusions

We proposed an MFG to model sustainable tourism based on the Jacobi-type processes. The MFG was exactly solvable under certain conditions where the controlled process became a mean field Jacobi process whose drift coefficient depended on the mean value of the process. Sustainable tourism and overtourism states were classified according to the boundary behavior of the derived controlled process. A finite difference method was applied to the cases in which no closed-form solution was found. The implications of the modeling results for Shiramine Village in Japan were discussed.

The proposed model is conceptual and should be extended to more realistic ones in the future. For example, tourism dynamics depend on the seasonality of tourist spots, which motivates us to consider a control problem with temporal periodicity. Another issue to be addressed is the coupling of a Jacobi-type process with the local natural resources and environmental dynamics that contribute to the local tourism industry. In this case, a higher-dimensional control problem must be addressed, whose closed-form solution would be more difficult to find than that of the proposed model. Therefore, a computational method will become essential. For example, a numerical solver, such as a machine-learning scheme would then become essential. Field surveys at the study site will continue to collect additional information and data to establish advanced mathematical models with quantified parameters and coefficients.

**Appendix**

*Proof of Proposition 1*

By definition, the right-hand side of (34) solves the HJB equation (28). The corresponding $u^*, w^*$ is given by (37) according to (26) and (27). The controlled process in (35) then becomes a classical Jacobi process with the source (due to **Lemma 1** and the assumption that $m^* \in D$ solves (33))

$$2M(m^*)E^* > 0 \tag{83}$$

and mean reversion

$$b - 2M(m^*)A^* > 0. \tag{84}$$

Moreover, due to $m^* \in D$ and **Lemma 1**, we have

$$b - 2M(m^*)A^* > 2M(m^*)E^*. \tag{85}$$

The inequalities (83)–(85) show that the controlled process $X$ is contained in $\bar{D}$ a.s. for $t > 0$ with the stationary average $m^*$ such that (33). The admissibility of $u^*$ and $w^*$ follows from their functional forms and the following inequalities:

$$\begin{aligned}
\mathbb{E}_{\mathbb{Q}}\left[\int_0^{+\infty} e^{-\delta s}(u_s^*)^2 \, ds\right] &= \mathbb{E}_{\mathbb{Q}}\left[\int_0^{+\infty} e^{-\delta s}\left(\frac{1}{c(m^*)}(A^* X_s + E^*)\right)^2 ds\right] \\
&\leq \left(\frac{1}{c(m^*)}(|A^*| + E^*)\right)^2 \mathbb{E}_{\mathbb{Q}}\left[\int_0^{+\infty} e^{-\delta s} ds\right] \\
&\leq \frac{1}{\delta}\left(\frac{1}{c(m^*)}(|A^*| + E^*)\right)^2 \\
&< +\infty
\end{aligned} \tag{86}$$

and

$$\begin{aligned}
\mathbb{E}_{\mathbb{Q}}\left[\int_0^{T+\infty} e^{-\delta s} g(X_s)(w_s^*)^2 \, ds\right] &= \mathbb{E}_{\mathbb{Q}}\left[\int_0^{+\infty} e^{-\delta s} X_s(1-X_s)\left(-\frac{\sigma}{\eta}\frac{1}{\sqrt{X_s(1-X_s)}}(A^* X_s + E^*)\right)^2 ds\right] \\
&\leq \left(\frac{\sigma}{\eta}(|A^*| + E^*)\right)^2 \mathbb{E}_{\mathbb{Q}}\left[\int_0^{+\infty} e^{-\delta s} ds\right] \\
&\leq \frac{1}{\delta}\left(\frac{1}{c(m^*)}(|A^*| + E^*)\right)^2 \\
&< +\infty
\end{aligned} \tag{87}$$

Based on the above discussion, it remains to be shown whether (34) yields the value function. This follows from a standard verification argument (e.g., Proof of Proposition 1 in Miao and Rivera, 2016; Proof of Theorem 1 in Lin and Riedel, 2021; Appendix in Pourbabaee, 2022) based on the twice continuous differentiability of the right-hand side of (34), inequalities (86) and (87), and another boundedness of the utility term given below:

$$\mathbb{E}_{\mathbb{Q}}\left[\int_0^{+\infty} e^{-\delta s}\left|f\left(X_s\right)\right| \mathrm{d}s\right] \leq \frac{1}{\delta}\max_{x\in\bar{D}}\left\{\left|-\frac{P}{2}x^2+Qx\right|\right\}<+\infty. \tag{88}$$

□

**Discretization of the finite difference method**

We present the finite difference method to discretize system (78)–(81). The time points are set as $\tau_i = i\Delta\tau$ with a constant $\Delta\tau > 0$. The space points are set as $x_j = j\Delta x$ ($j = 0,1,2,...,I_x$) and $\Delta x = 1/I_x$, with some natural numbers $I_x \geq 2$. The discretized $\Phi$ and $p$ at the point $(\tau_i, x_j)$ are expressed as $\Phi_{i,j}$ and $p_{i,j}$, respectively. For simplicity, we chose the initial guesses $\Phi_{0,j} = 0$ and $p_{0,j} = 1$ ($j = 0,1,2,...,I_x$). The discretization scheme in our setting is an adaptation of that proposed by Achdou and Capuzzo-Dolcetta (2010). A difference between their discretization schemes and ours is that the former deals with a system on a torus that does not have any boundaries, whereas the latter has left and right boundaries. For brevity, we set $\Phi_{i,-1} = \Phi_{i,I_x+1} = 0$.

The HJB equation in the systems (78)–(81) is discretized as follows:

$$\frac{\Phi_{i+1,j}-\Phi_{i,j}}{\Delta\tau} = -\delta\Phi_{i,j} + \mathrm{Ham}_{i,j} - bx_j\frac{\Phi_{i,j}-\Phi_{i,j-1}}{\Delta x} \\ + \frac{\sigma^2}{2}x_j(1-x_j)\frac{\Phi_{i,j-1}-2\Phi_{i,j}+\Phi_{i,j+1}}{(\Delta x)^2} + f(x_j), \quad i=0,1,2,... \text{ and } j=0,1,2,...,I_x \tag{89}$$

with the numerical Hamiltonian $\mathrm{Ham}_{i,j}$ evaluated as

$$\mathrm{Ham}_{i,j} = M(m_i) \times \begin{cases} \left(\dfrac{\Phi_{i,1}-\Phi_{i,0}}{\Delta x}\right)^2 & (j=0) \\ \left(\max\left\{\dfrac{\Phi_{i,j+1}-\Phi_{i,j}}{\Delta x},0\right\}\right)^2 + \left(-\min\left\{\dfrac{\Phi_{i,j}-\Phi_{i,j-1}}{\Delta x},0\right\}\right)^2 & (1\leq j\leq I_x-1) \\ \left(\dfrac{\Phi_{i,I_x}-\Phi_{i,I_x-1}}{\Delta x}\right)^2 & (j=I_x) \end{cases} \tag{90}$$

Here, the discretized $m$ is given by

$$m_i = \Delta x\left(\frac{1}{2}x_0 p_{i,0} + \sum_{j=1}^{I_x-1}x_j p_{i,j} + \frac{1}{2}x_{I_x}p_{i,I_x}\right), \quad i=0,1,2,.... \tag{91}$$

The factor $\dfrac{1}{2}$ in (91) implies that the control volumes at the boundaries have half the area of those of the interior. The FP equation is discretized as follows:

$$\frac{p_{i+1,j}-p_{i,j}}{\Delta\tau} = -\frac{\mathrm{Flux}_{i,j}-\mathrm{Flux}_{i,j-1}}{\Delta x}, \quad i=0,1,2,..., \quad j=1,2,3,...,I_x \tag{92}$$

with the numerical flux $\mathrm{Flux}_{i,j}$ ($\mathrm{Flux}_{i,0} = \mathrm{Flux}_{i,I_x} = 0$) given by

$$\text{Flux}_{i,j} = 2\left( m_i \min\left\{ \frac{\Phi_{i,j+1} - \Phi_{i,j}}{\Delta x}, 0 \right\} p_{i,j+1} + m_i \max\left\{ \frac{\Phi_{i,j+1} - \Phi_{i,j}}{\Delta x}, 0 \right\} p_{i,j} \right)$$
$$-bx_j p_{i,j+1} - \frac{\sigma^2}{2} \frac{x_{j+1}(1-x_{j+1}) p_{i,j+1} - x_j(1-x_j) p_{i,j}}{\Delta x}, \quad j=1,2,3,...,I_x-1. \quad (93)$$

The discretized system (89)-(90) is computed from $i = 0,1,2,...$, starting with a prescribed initial guess. A stationary solution is obtained computationally at time $\tau_i$ if the following condition is satisfied for a small error threshold $\varepsilon > 0$:

$$\max_{j=0,1,2,...,I_x} \left\{ \left| \Phi_{i,j} - \Phi_{i-1,j} \right|, \left| p_{i,j} - p_{i-1,j} \right| \right\} \leq \varepsilon. \quad (94)$$

We chose $\varepsilon = 10^{-12}$. For the ergodic control case ($\delta = 0$), we use the convergence criterion

$$\max_{j=0,1,2,...,I_x} \left\{ \left| p_{i,j} - p_{i-1,j} \right| \right\} \leq \varepsilon. \quad (95)$$

because in this case, we asymptotically have $\Phi_{i,j} - \Phi_{i-1,j} = \Delta\tau H$ for a sufficiently large $i$ when an effective Hamiltonian $H$ exists. In this case, numerically verifying the stationarity of the PDF is sufficient.

For small values of $\Delta\tau > 0$, the computed PDF remains non-negative and satisfies the probability conservation law

$$1 = \Delta x \left( \frac{1}{2} p_{i,0} + \sum_{j=1}^{I_x-1} p_{i,j} + \frac{1}{2} p_{i,I_x} \right), \quad i = 0,1,2,.... \quad (96)$$

The proof is omitted here because it essentially follows from the stability argument for time-explicit discretization methods (e.g., Neena et al, 2024; Richter et al., 2023), where we assume that the differences $\frac{\Phi_{i,j+1} - \Phi_{i,j}}{\Delta x}$ and $\frac{\Phi_{i,j} - \Phi_{i,j-1}}{\Delta x}$ are strictly bounded in both space and time. In addition, owing to the use of forward Euler discretization in time, we need to choose a sufficiently small $\Delta\tau$ value such that $\Delta\tau = O\left((\Delta x)^2\right)$ ensures the non-negativity of the PDF $p$. One may employ implicit discretization in time to improve the computational efficiency of the finite difference method; however, we do not address this issue because we are interested in stationary solutions rather than transient ones.

**Convergence of the finite difference method**

The proposed finite difference method uses monotone discretization for both the HJB and FP equations; hence, it has first-order accuracy (at most) with respect to discrete norms (Almulla et al., 2017; Osborne and Smears, 2024). Therefore, the convergence study focuses on whether the finite-difference scheme can reasonably approximate the value function $\Phi$, PDF $p$, and average $m$. We examine cases using $\delta = 0.1$ and $\delta = 0$. The computed $\Phi$, $p$, and $m$ values are compared with the closed-form solutions derived in **Proposition 1**. The other parameter values are the same as those in **Section 4**. We fix $\Delta\tau = 0.00025$, which is sufficiently small to compute the numerical solutions stably.

As shown in **Tables A1-A2**, the convergence of numerical solutions to closed-form solutions is first order in $I_x$, and the errors are sufficiently small at $I_x = 400$. The errors between the computed and closed-form $m$ are sufficiently small at $I_x = 400$; e.g., for $\delta = 0.1$, $m = 0.48048$ for the closed-form solution, and $m = 0.47999$ for the numerical solution, with a relative error of 0.1%. Therefore, we used the resolution $I_x = 400$ in the main text. Numerical solutions are slightly diffused compared with the closed-form solutions owing to numerical diffusion and do not contain spurious oscillations owing to the use of monotone discretization (see **Figures A1-A4**).

**Table A1.** $l^1$ errors between the closed-form and computed $\Phi$ and $p$: $\delta = 0.1$ and $\eta = 1$.

| $I_x$ | $\Phi$ | $p$ |
|---|---|---|
| 50 | 0.005754 | 0.08649 |
| 100 | 0.002699 | 0.04559 |
| 200 | 0.001282 | 0.02345 |
| 400 | 0.000617 | 0.01190 |

**Table A2.** $l^1$ errors between the closed-form and computed $\Phi$ and $p$: $\delta = 0$ and $\eta = 0.1$. We normalized $\Phi$ such that $\Phi(0) = 0$.

| $I_x$ | $\Phi$ | $p$ |
|---|---|---|
| 50 | 0.004279 | 0.07828 |
| 100 | 0.002513 | 0.03955 |
| 200 | 0.001403 | 0.01990 |
| 400 | 0.000761 | 0.00999 |

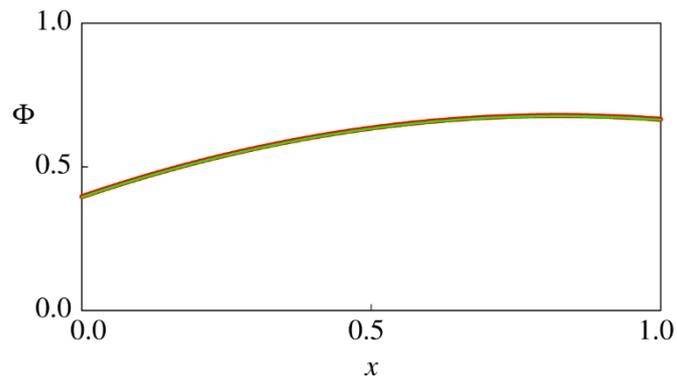

**Fig. A1** Comparison of the closed-form and computed $\Phi$ plots for $\delta = 0.1$ and $\eta = 1$ (closed-form solution (red) and numerical solution with $I_x = 200$ (green))

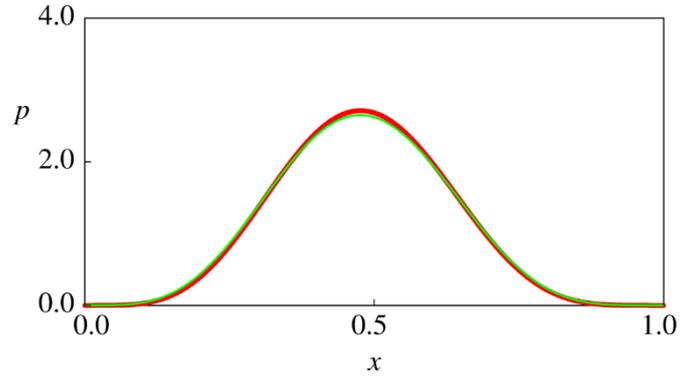

**Fig. A2** Comparison of the closed-form and computed $p$ plots for $\delta = 0.1$ and $\eta = 1$ (closed-form solution (red) and numerical solution with $I_x = 200$ (green))

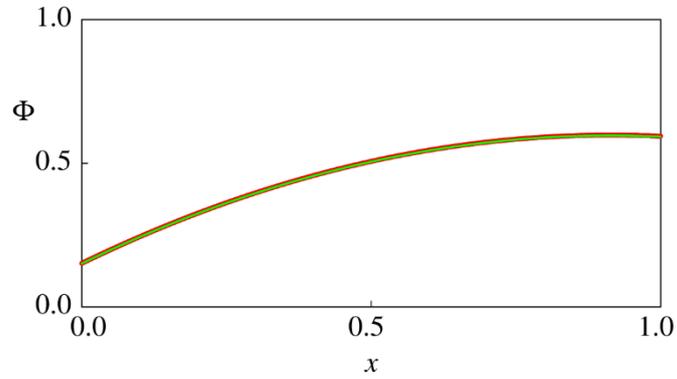

**Fig. A3** Comparison of the closed-form and computed $\Phi$ plots for $\delta = 0.1$ and $\eta = 0.1$ (closed-form solution (red) and numerical solution with $I_x = 200$ (green))

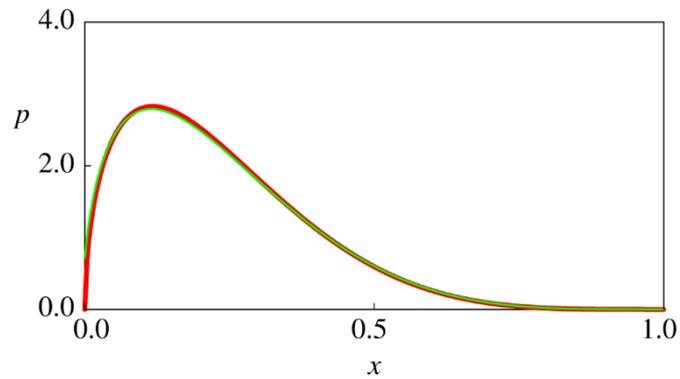

**Fig. A4** Comparison of the closed-form and computed $p$ plots for $\delta = 0.1$ and $\eta = 0.1$ (closed-form solution (red) and numerical solution with $I_x = 200$ (green))

Sass, G. G., Feiner, Z. S., & Shaw, S. L. (2021). Empirical evidence for depensation in freshwater fisheries. Fisheries Magazine, 46(6), 266-276. https://doi.org/10.1002/fsh.10584

Sun, M., & Wang, C. (2024). Asymptotic behavior of Riemann solutions for the one-dimensional mean-field games in conservative form with the logarithmic coupling term. International Journal of Non-Linear Mechanics, 104837. https://doi.org/10.1016/j.ijnonlinmec.2024.104837

Takahashi, T. (2024). The conflict between residents and tourists: on the variety-shifting effect of tourism growth. The Japanese Economic Review, 75(1), 121-145. https://doi.org/10.1007/s42973-021-00108-5

Tong, Z., & Liu, A. (2022). Pricing variance swaps under subordinated Jacobi stochastic volatility models. Physica A: Statistical Mechanics and its Applications, 593, 126941. https://doi.org/10.1016/j.physa.2022.126941

Vanelli, F. M., & Kobiyama, M. (2021). How can socio-hydrology contribute to natural disaster risk reduction?. Hydrological Sciences Journal, 66(12), 1758-1766. https://doi.org/10.1080/02626667.2021.1967356

Xepapadeas, A. (2024). Uncertainty and climate change: The IPCC approach vs decision theory. Journal of Behavioral and Experimental Economics, 109, 102188. https://doi.org/10.1016/j.socec.2024.102188